\documentclass[10pt,twoside]{article}
\usepackage{graphicx}
\usepackage{amsmath,amssymb,amscd,enumerate}
\usepackage{Latex-document}

\newcommand{\eps}{\varepsilon}
\newcommand{\diff}[2]{\frac{\partial #1}{\partial #2}}
\def\volumeno{I}
\title{\bf  Energy Landscapes and
Rare Events \vskip 6mm}
\author{{\bf Weinan E}\thanks{Department of Mathematics and
PACM, Princeton University, Fine Hall, Princeton, NJ 08544, USA
and School of Mathematics, Peking University, Beijing, 100871,
China. E-mail: weinan@math.princeton.edu} \quad Weiqing
Ren\thanks{Courant Institute of Mathematical Sciences, New York
University, New York 10012, USA. Email: weiqing@cims.nyu.edu}
\quad Eric Vanden-Eijnden\thanks{Courant Institute of Mathematical
Sciences, New York University, New York 10012, USA. Email:
eve2@cims.nyu.edu}\vspace*{-0.5cm}}
\date{\vspace{-8mm}}

\begin{document}
\maketitle

\thispagestyle{first} \setcounter{page}{621}

\begin{abstract}\vskip 3mm
Many problems in physics, material sciences, chemistry and
biology can be abstractly formulated as a system that navigates
over a complex energy landscape of high or infinite dimensions.
Well-known examples include phase transitions of condensed matter,
conformational changes of biopolymers, and chemical reactions.
The energy landscape typically exhibits multiscale features,
giving rise to the multiscale nature of the dynamics. This is one
of the main challenges that we face in computational science.
In this report, we will review the recent work done by scientists from
several disciplines on probing such energy landscapes. Of particular
interest is the  analysis and computation of
transition pathways and transition rates between metastable states.
We will then present the string method that has proven to be very
effective for some truly complex systems in material science and
chemistry.

\vskip 4.5mm


\end{abstract}

\vskip 12mm

\section{Introduction} \label{section 1}\setzero
\vskip-5mm \hspace{5mm }

Many problems in biology, chemistry and material science can be
formulated as the study of the energy or free energy landscape of the
underlying system. Well-known examples of such problems include the
conformational changes of macromolecules, chemical reactions and
nucleation in condensed systems. Very often the dimension of the state
space is very large, and the energy landscape exhibits a hierarchy of
structures and scales.
These problems are becoming a
major challenge in their respective scientific disciplines and are
beginning to receive attention from the mathematics community. In this
article, we report recent work
in this direction. For a detailed account, we refer to
\cite{ERV1,ERV2,ERV3,ERV4,Ren}.

We begin with a simple example. Plotted in Figure 1 is the solution of
the stochastic differential equation
\begin{equation}
\label{eq:1}
dx(t)=-\nabla_x V(x(t))dt+\sqrt{\eps}dW(t)
\end{equation}
where the potential
\begin{equation}
\label{eq:2}
V(x)=\frac14(1-x^2)^2
\end{equation}
and $dW(t)$ is Gaussian white noise, $\eps=0.06,x(0)=-1$. Without the
random perturbation, the solution would be $x(t)\equiv x(0)=-1$. Indeed
the deterministic part of dynamics in \eqref{eq:1} does nothing but
taking the system to local equilibrium states. With the random
perturbation, the solution, over long time, exhibits completely different
behavior. It fluctuates around the two local minima of $V,x=-1$ and 1,
with sudden transitions between these two states. The time scale of the
transition, $t_M$ is much larger than the time scale of the fluctuation
around the local minima, $t_R$. For this reason, we refer to $x=-1$ and
1 as the metastable states.

\begin{center}
\includegraphics[width=4.5in,height=3in]{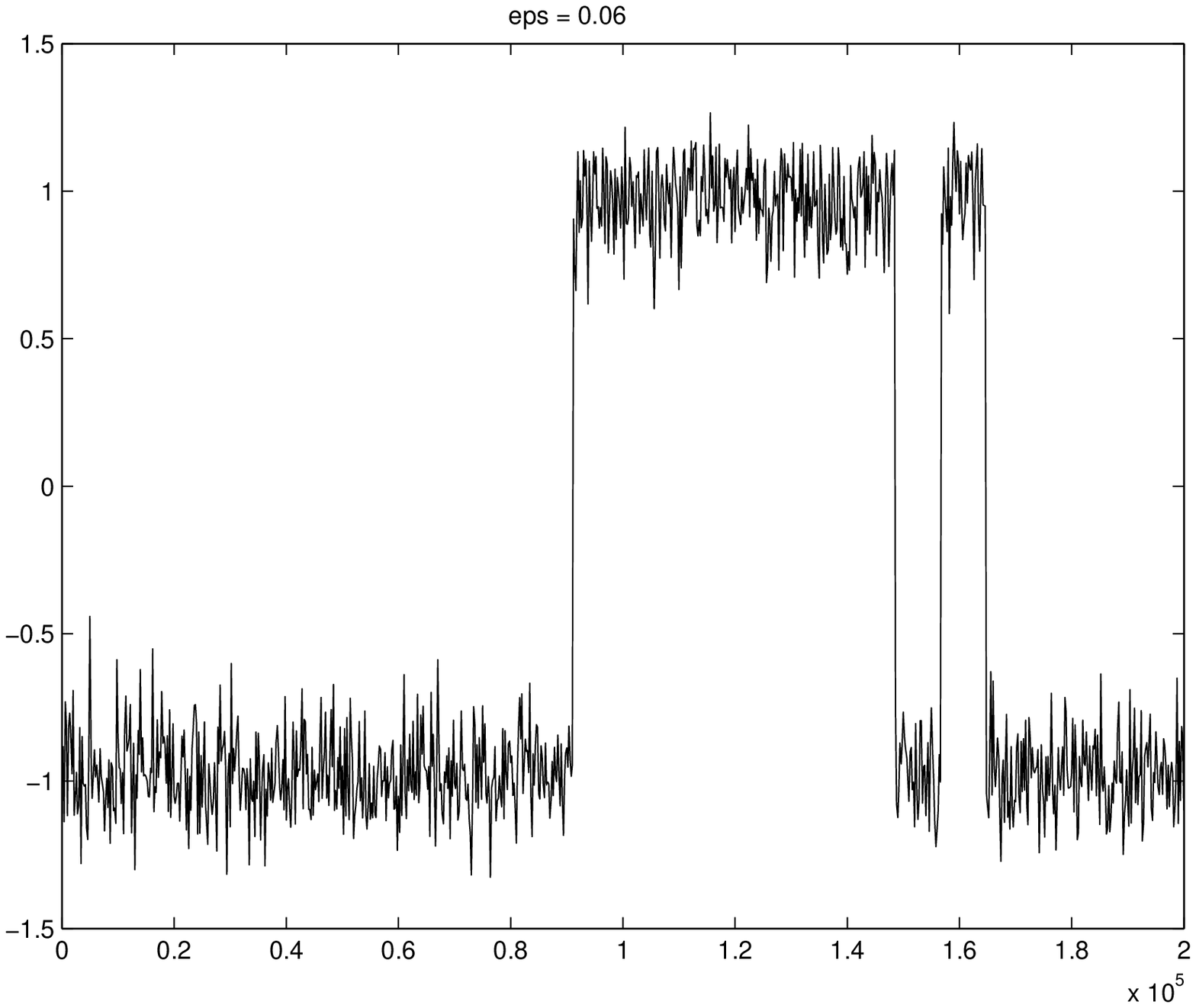}
\begin{minipage}[h]{10cm}
Figure 1. Time series of the solution to the stochastic
differential equation (1.1), with $\eps=0.06$.
\end{minipage}
\end{center}
\vspace{.3in}

Obviously the transition between the metastable states is of more
interest than the local fluctuation around them. The
transition time is much larger since it requires the system to overcome
the energy barrier between the two states. This is only possible because
of the noise. When $\eps$ is small, a huge noise term is required to
accomplish this. For this reason, such events are very rare, and this is
the origin of the disparity between the time scales $t_M$ and $t_R$.

This simple example illustrates one of the major difficulties in modeling
such systems, namely the disparity of the time scales. It does not,
however, illustrate the other major difficulty, namely, the large
dimension of the state space and the complexity of the energy
landscapes. Indeed for typical systems of interests the energy landscape
can be very complex. There can be a huge number of
local minima in the state space. The usual concept of hopping over
barriers via saddle points may not apply (see \cite{deboge01}).

In applications, the noise comes typically from thermal noise. In this
case, we should note that even though the potential energy landscapes
might be rough and contain small scale features, the system itself
experiences a much smoother landscape, the free energy landscape, since
some of the small scale features on the potential energy landscape
 are smoothed out by the thermal noise.

Our objective in modeling such systems are the following:

\begin{enumerate}
\item Find the transition mechanism between the metastable states.
\item Find the transition rates.
\item Reduce the original dynamics to the dynamics of a Markov chain on
the metastable states.
\end{enumerate}

Our discussion will be centered around the following model problems:
\begin{equation}
\label{eq:3}
\gamma\dot{x}(t)=-\nabla V(x(t))+\sqrt{\eps}\dot{W}(t)
\end{equation}
or
\begin{equation}
\label{eq:4}
m\ddot{x}(t)+\gamma\dot{x}(t)=-\nabla V(x(t))+\sqrt{m\eps}\dot{W}(t)
\end{equation}
$\eps$ is related to the temperature of the system by $ \eps
=2 \gamma k_B T$ where $k_B$ is the Boltzmann constant.
We refer to \eqref{eq:3} as type-I gradient flow and \eqref{eq:4} as
type-II gradient flow.

Before proceeding further, let us remark that there is a very
well-developed theory, the large-deviation theory, or the
Wentzell-Freidlin theory \cite{FW}, that deals precisely with
questions of the type that we discussed above. However as was explained
in \cite{ERV4,Ren}, this theory is not best suited for numerical purpose.
Therefore we will seek an alternative theoretical framework that is
more useful for numerical computations.

\section{Transition state theory} \label{section 2}
\setzero\vskip-5mm \hspace{5mm }

Transition state theory (TST) \cite{hatabo90}
has been the classical framework for
addressing the questions we are interested in. It assumes the existence
and explicit knowledge of a reaction coordinate, denoted by $q$, that
connects the two metastable states. In addition it assumes that along
the reaction coordinate there exists a well-defined transition state,
which is typically the saddle point configuration, say at $q=0$, and the
two regions $\{q<0\}$ and $\{q>0\}$ defines the two metastable regions
$A$ and $B$. For these reasons, transition state theory is restricted to
cases when the system is simple and the energy landscape is smooth,
i.e. the energy barriers are larger than the thermal energy $k_B T$.

Knowing the transition state, TST calculates the transition rates by
placing particles at the transition state, and measuring the flux that
goes into the two regions. For example, the transition rate from $A$ to
$B$ is given approximately by
\begin{equation}
\label{eq:5}
k_{A\to
B}=\frac1{Z_0}\int\dot{q}(t)\delta_\Gamma(q(t))\theta(q(t))d\mu_A(q(0))
\end{equation}
where
\begin{equation}
\label{eq:6}
Z_0=\int d\mu_A(q(0))
\end{equation}
Here $\delta_\Gamma$ is the surface delta function at $q=0$, $\theta$ is
the Heaviside function, $\mu_A$ is the Lebesgue measure restricted to
$A$. For a system with a single particle of mass $m$ and potential
$V$, this gives \cite{hatabo90}
\begin{equation}
\label{eq:7}
k_{A\to B}=\frac{\omega_0}{2\pi}e^{-\frac{\delta E}{k_BT}}
\end{equation}
where $\delta E$ is the energy barrier at the transition state,
$\omega_0=\left(\frac{V''(x_A)}m\right)^{\frac12}$, $x_A$ denotes the
location of the local minimum inside $A$. Formulas such as \eqref{eq:7}
are the origin of the Arrhenius law for chemical reaction rates and
Boltzmann factor for hopping rates in kinetic Monte Carlo models.

\section{Reduction to Markov chains on graphs} \label{section 3}
\setzero\vskip-5mm \hspace{5mm }

For simplicity, we will discuss mainly type-I gradient flows
\eqref{eq:3}. The Fokker-Planck equation can be expressed as
\begin{equation}
\label{eq:8}
\diff{p}{t}(x,t)=\nabla\cdot\left(p_s(x)\nabla\left(\frac{p(x,t)}{p_s(x)}\right)\right)
\end{equation}
where $p_s$ is the equilibrium distribution
$$p_s(x)=\frac1Ze^{-\frac{V(x)}{k_BT}}$$
$Z$ is the normalization constant
$Z=\int_{R^n}e^{-\frac{V(x)}{k_BT}}dx$. The states of the Markov chain
consist of the sets $\{B_j\}^J_{j=1}$, where $\{B_j\}^J_{j=1}$
satisfies
\begin{enumerate}
\item The $B_j$'s are mutually disjoint
\item
\begin{equation}
\label{eq:9}
\int_Bp_s(x)dx=1+o(k_BT)
\end{equation}
\end{enumerate}
where $B=\cup^J_{j=1}B_j$. An illustration of the collection the
$\{B_j\}^J_{j=1}$ is given in Figure 2. $\{B_j\}$ depends on $T$. As $T$
decreases, the $B_j$'s exhibits a hierarchical structure.

\begin{center}
\includegraphics[width=4in,height=2.5in]{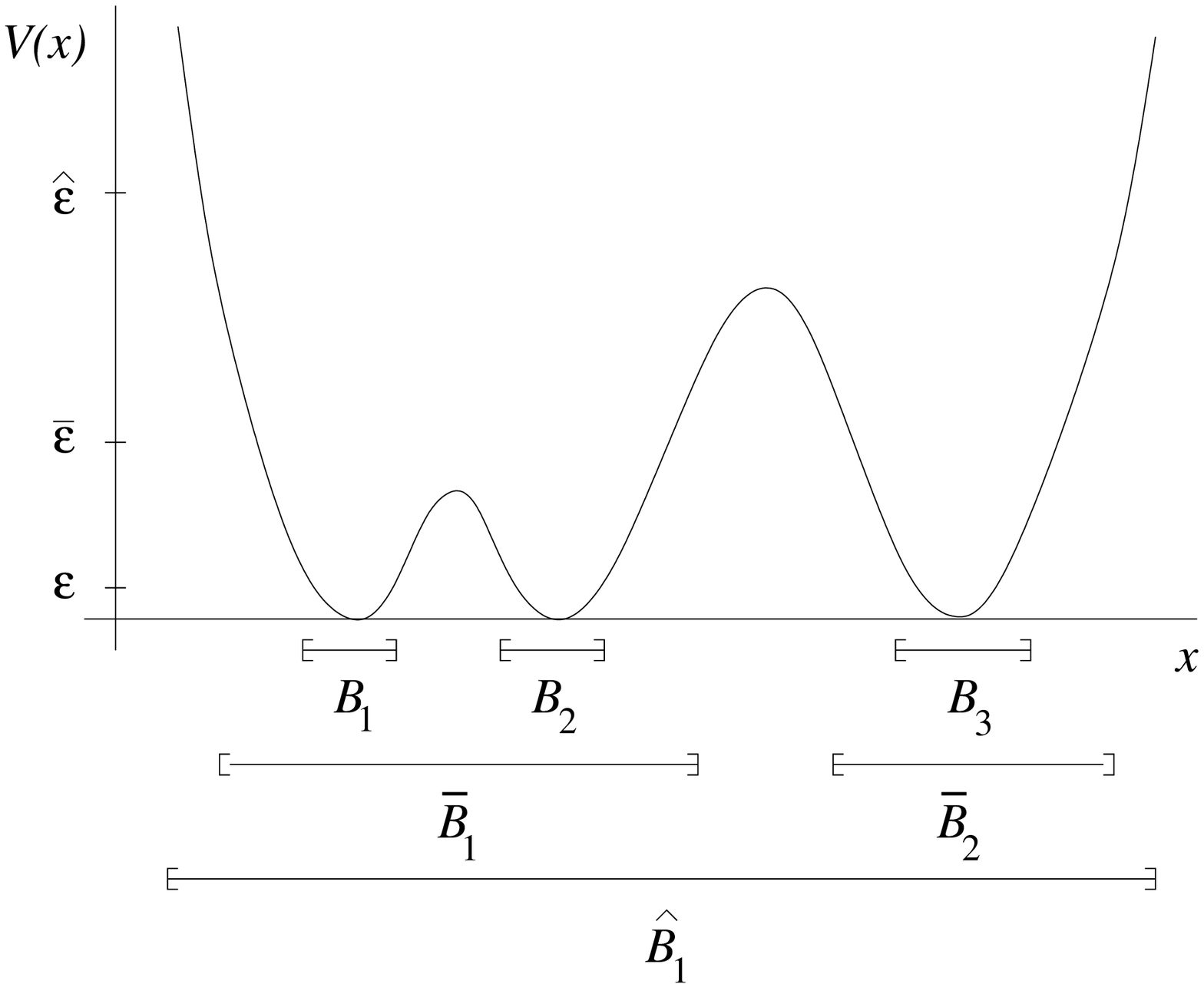}
\begin{minipage}[h]{10cm}
Figure 2. Illustration of the collection of metastable sets
$\{B_j\}$ at different energies.
\end{minipage}
\end{center}
\vspace{.2in}

Having defined the states of the Markov chain, we next compute the
transition rates between neighboring states.
Denote by $A$ and $B$ two such neighboring states. We would like to compute
the transition rate from $A$ to $B$. Without loss of generality, we may
assume $J=2$. Let $B_1,B_2$ be the metastable region containing $A$ and
$B$ respectively, and let $n_j(t)=\int_{B_j}p(x,t)dx$,
$N_j=\int_{B_j}p_s(x)dx$, $j=1,2$. Applying Laplace's
method to (3.1) we get \cite{ERV4}
\begin{equation}
\label{eq:11}
\dot{n}_j(t)=\frac\eps{\kappa}\left(\frac{n_2(t)}{N_2}-
\frac{n_1(t)}{N_1}\right)+\mbox{higher
order terms}
\end{equation}
where
\begin{equation}
\label{eq:12}
\kappa=\int^1_0d\alpha\left(\int_{S^0(\alpha)}p_s(x)dx\right)^{-1}|\varphi^0(\alpha)|
\end{equation}
$\{\varphi^0(\alpha),0\le\alpha\le1\}$ is a so-called minimal energy
path, to be defined below, $\{S^0(\alpha)\}$ is the family of
hyperplanes normal to $\varphi^0$.

The minimal energy path (MEP) is defined as follows. If $V$ is smooth,
then $\varphi^0$ is a MEP if
\begin{equation}
\label{eq:13}
(\nabla V)^\perp(\varphi^0(\alpha))=0
\end{equation}
for all $\alpha\in[0,1]$, i.e. $\nabla V$ restricted to $\varphi^0$ is
parallel to $\varphi^0$. In general there is not a unique $\varphi^0$
that is particularly significant, but rather a collection (a tube
or several tubes) of paths
contribute to the transition rates. However, one can define a MEP
self-consistently via the equation
\begin{equation}
\label{eq:14}
\varphi^0(\alpha)=\frac1{Z(\alpha)}\int_{S^0(\alpha)}xe^{-\frac{V(x)}{k_BT}}dx
\end{equation}
where $Z(\alpha)=\int_{S^0(\alpha)}e^{-\frac{V(x)}{k_BT}}dx$. In the
case when $V$ has two scales: $V=\bar{V}+\delta V$,  $|\delta
V|\le O(k_BT)$, and $\bar{V}$ is smooth, then $\varphi^0$ can be defined
as the MEP of $\bar{V}$.

For type-I gradient systems, if the two metastable sets are separated by a
single saddle point, then the MEP is the unstable manifold associated
with the saddle point. MEP for type-II gradient systems is less trivial.
In this case (3.3)-(3.6) have to be modified \cite{ERV4}.
Consider the simple example
$$\left\{\begin{array}{ll}
\dot{q} & =p \\
\dot{p} & =-\frac{\partial V}{\partial q}(q)- p+\sqrt\eps\dot{W}
\end{array}\right.$$
where $V(q)=\frac14 (1-q^2)^2$.
It has two local equilibrium states $A=(-1,0)$ and $B=(1,0)$.
The MEP that connects
$A$ and $B$ is plotted in Figure 3. It is not a smooth curve. The
velocity is reversed at the saddle point. To verify that the MEP does
reflect the true behavior of the transition path, we also plot the
transition path obtained from direct simulation of the stochastic
differential equation.

\begin{center}
\includegraphics[width=4.5in,height=3in]{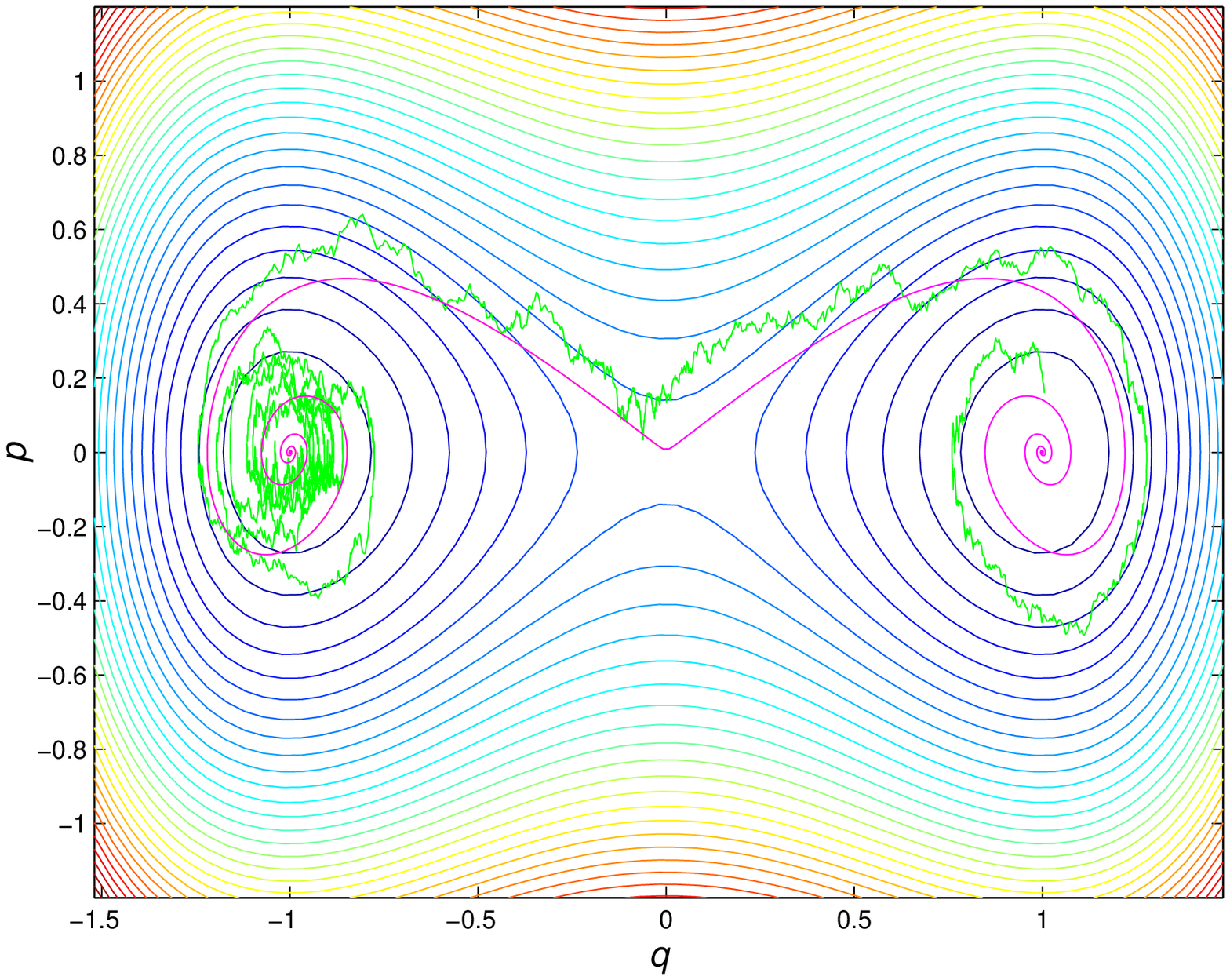}
\begin{minipage}[h]{10cm}
Figure 3. MEP for type-II gradient systems. The  red line is the
MEP in phase space, the green  line is the transition path
computed from solving the stochastic differential equation.
\end{minipage}
\end{center}
\vspace{.2in}

MEP is a very important concept since it defines the ``most probable''
transition path from which transition rates can be computed via equation
\eqref{eq:11}. However we should emphasize that from a numerical point
of view our task is not that of a conventional optimization or control
problem, since there is not an objective function that we can easily
work with. Instead our aim is to perform importance sampling in path space
to sample the paths that contribute significantly to the switching.

Finally, if there exists a MEP that connects two metastable sets  without
going through a third one, then we connect these two metastable sets
by a link. In this way, we form a graph. The original dynamics is then
reduced to a Markov chain on this graph.

\section{Previous numerical techniques} \label{section 4}
\setzero\vskip-5mm \hspace{5mm }

A variety of numerical techniques have been developed, most
prominently in chemistry, but also in biology and material
science, for computing MEPs and sometimes transition rates. Among
the most well-known techniques in the chemistry literature are the
nudged elastic band method and the transition path sampling
technique. The nudged elastic band method (NEB) \cite{jomijo98}
aims at computing the MEP defined by (3.5). It represents the MEP
by a discrete chain of states. These states evolve according to
the potential forces of the system. To prevent the states from all
falling to the two local equilibrium states, a spring force is
applied to neighboring states to penalize the non-uniformity in
the distribution of the states along the chain. This by itself may
cause convergence to a path which is not a MEP. Hence a nudging
technique is used, namely only the normal component of the
potential force and the tangential component of the spring force
is applied. NEB is a very effective method for small systems with
relatively smooth energy landscapes. It has two main drawbacks.
One is that it is highly inefficient and may not even be
applicable to systems with rough energy landscapes. The other is
the choice of the elastic constant. A large elastic constant
requires a small time step in the evolution of the states. A small
elastic constant will not achieve the desired uniformity of the
states and hence will not give  the required accuracy for the
energy barrier.

A second important technique is transition path sampling (TPS)
\cite{debocs98,deboge01}. This
method aims at complex systems with rough energy landscapes
by developing a Monte Carlo technique that samples the path space. Its
efficiency hinges on the ability to produce new accepted paths from old
ones.

Other techniques include the ridge method, blue mooth sampling,
etc. \cite{cacihyka89}. Often these methods require knowing
beforehand the reaction coordinate.

Elber et. al propose to minimize the Onsager-Machlup action as a way of
finding the most probable path for macromolecular systems \cite{olel96}. The
Onsager-Machlup action is the same as the Wentzell-Freidlin action. From
a numerical point of view, there are certain difficulties associated
with minimizing this action functional. These issues are discussed in
\cite{ERV4}.

\section{The string method} \label{section 5}
\setzero\vskip-5mm \hspace{5mm }

The basic idea in the string method, developed in
\cite{ERV1,ERV2,ERV4}, is to represent transition paths by their
intrinsic parameterization in order to efficiently evolve and
sample paths in path space. It has two versions. The zero
temperature version is designed for smooth energy landscapes. The
finite temperature version is designed for rough energy landscapes
in which case thermal noise acts to smooth out the small scale
features.

The simplest example of a zero-temperature string method is to evolve
curves in path space by the gradient flow
\begin{equation}
\label{eq:5.11}
\varphi_t(\alpha,t)=-(\nabla
V)^\perp(\varphi(\alpha,t))+r(\alpha,t)\hat{\tau}(\alpha,t)
\end{equation}
Here $\hat{\tau}$ is the tangent vector of the curve
$\{\varphi(\cdot,t)\}$, $(\nabla V)^\perp(\varphi)$ denotes the
component of $\nabla V$ normal to $\hat{\tau}$, $r$ is the Lagrange
multiplier that enforces certain specific parameterization of the curves.
For example if we require equal arclength parameterization, then we need
$\diff{}{\alpha}|\varphi_\alpha|=0$, i.e.
\begin{equation}
\label{eq:5.12}
r(\alpha,t)  =  \alpha\int^1_0\nabla
V(\varphi(\alpha',t))\cdot\hat{\tau}(\alpha',t)d\alpha'
  -\int^\alpha_0\nabla
V(\varphi(\alpha',t))\cdot\hat{\tau}(\alpha',t)d\alpha'
\end{equation}
We call such curves with intrinsic parameterization strings.

In practice the strings are discretized into a collection of points.
These points move according to the normal component of the potential
force. A reparameterization step is applied once in a while to enforce
the proper parameterization of the strings.

The finite temperature string method is designed for systems with rough
energy landscapes, particularly the case when the potential can be
expressed in the form
$$V(x)=\bar{V}(x)+\delta V(x)$$
where $\bar{V}$ is smooth and $|\delta V|\le O(k_BT)$. In this case we
would like to compute the MEP of $\bar{V}$ without first computing
$\bar{V}$ explicitly. This is achieved by creating an ensemble of a
special type. Our computational object will be a string connecting the
two metastable sets,
together with a family of probability measures on the
hyperplanes normal to the string. Consider the stochastic equation
\begin{equation}
\label{eq:5.1}
\varphi^\omega_t(\alpha,t)=-{\bf P}^0_\alpha(\nabla
V(\varphi^\omega(\alpha,t)))+r^0(\alpha,t)\hat{\tau}^0(\alpha,t)+
{\bf P}^0_\alpha\eta^\omega(\alpha,t)
\end{equation}
where $\eta^\omega$ is Gaussian noise with mean 0 and correlation
$${\bf E}\eta^\omega(\alpha,t)\eta^\omega(\alpha',\tau)=\left\{\begin{array}{ll}
2k_BT\delta(t-\tau), & \mbox{if }\alpha=\alpha' \\
0, & \mbox{if }\alpha\not=\alpha'\end{array}\right.$$
The projection operator ${\bf P}^0_\alpha$ is defined by projecting to the
hyperplane normal to the string
$\{\varphi^0(\alpha,t),0\le\alpha\le1\}$, where
$$\varphi^0(\alpha,t)={\bf E}\varphi^\omega(\alpha,t)$$
$\hat{\tau}^0(\alpha,t)$ is the tangent vector of $\varphi^0$ at
$\alpha$, $r^0$ is the Lagrange multiplier that enforces proper
parameterization of $\varphi^0$.

{\bf Theorem 1.1.} \it
\begin{enumerate}
\item The statistical steady state of \eqref{eq:5.1} satisfies:
\begin{equation}
\label{eq:5.13}
\varphi^0(\alpha)=\frac1{Z(\alpha)}\int_{S^0(\alpha)}xe^{-\frac{V(x)}{k_BT}}dx
\end{equation}
where $S^0(\alpha)$ is the hyperplane normal to $\varphi^0$ at $\alpha$,
$Z(\alpha)=\int_{S^0(\alpha)}e^{-\frac{V(x)}{k_BT}}dx$.

\item The stationary distribution of \eqref{eq:5.1} is given by the family of
distributions
\begin{equation}
\label{eq:5.14}
\mu_\alpha(x)=\frac1{Z(\alpha)}e^{-\frac{V(x)}{k_BT}}\delta_{S^0(\alpha)}(x)
\end{equation}
\end{enumerate}
\rm

Knowing $\{\varphi^0\}$ and $\{\mu_\alpha\}$, the transition rates and
free energy landscapes can be computed, see \cite{ERV4}.

The finite temperature string method is applied to the perturbed Mueller
potential
$$V(x)=V_m(x)+\delta V(x)$$
where $V_m(x)$ is the so-called Mueller potential
(see \cite{olel96}), $\delta V$ is a random perturbation.
The results of the (finite
temperature) string method is shown in Figure 4. Also plotted is the
MEP of $V_m$ (since $\bar{V}=V_m$ is explicitly known for this
particular example) as well as the fluctuations around $\varphi^0$.

\begin{center}
\includegraphics[width=4.5in,height=2.5in]{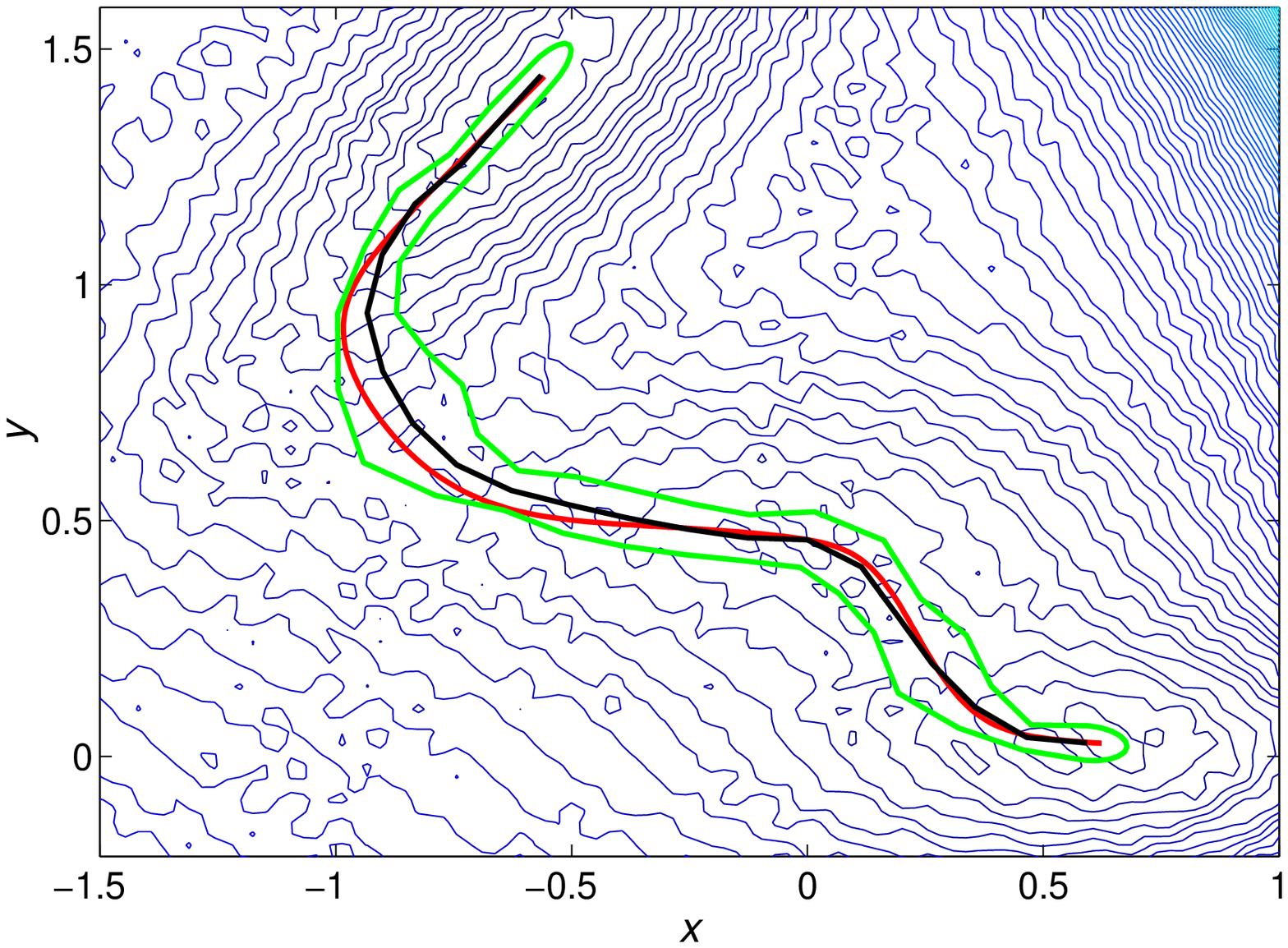}
\begin{minipage}[h]{10cm}
Figure 4. Effective MEP and local fluctuations for the perturbed
mueller potential. The red curve is the MEP for $\bar{V} = V_m$.
The black curve is the MEP computed from the finite temperature
string method. The green curves show the size of the fluctuations.
\end{minipage}
\end{center}
\vspace{.2in}

\section{Concluding remarks} \label{section 6}
\setzero\vskip-5mm \hspace{5mm }

There are several important topics that we did not cover in this brief
report. These include the effect of dynamics, non-gradient systems, and
acceleration techniques. These are discussed in
\cite{ERV1,ERV2,ERV4,Ren}. Also found in these references are
applications of the ideas discussed here to thermal activated reversal
of magnetic thin films, models of martensitic transformations, and the
formation of $C_{60}$ from 60 carbon atoms. The last example is a case
when the barrier is entropic. Even though the potential energy is mainly
going downhill, the free energy has barriers because of entropic effects.
Such examples are found frequently in biopolymers.

From a numerical point of view, our main idea for overcoming the
difficulty caused by the disparity of the times scales is to reformulate
the problem as a boundary value problem instead of initial  value
problem, since we have some knowledge of the initial and final
state of the system. Compared with other existing methods that
assume explicit knowledge of a reaction coordinate, our method
finds the reaction coordinate self-consistently during the computation.

The topic discussed here is relatively new in applied mathematics, but
it is of paramount importance in science and particularly computational
science. Progress in this area will likely have a fundamental impact in
many areas of applications.

\vspace{.15in}\noindent
{\bf Acknowledgment.}
We are  grateful to Bob Kohn for many stimulating
discussions. This work is partially supported by NSF via grant
DMS01-30107.

\label{lastpage}


\begin{thebibliography}{aa}


\bibitem{cacihyka89} Carter, E. A.; Cicotti, G.; Hynes, J. T.; Kapral,
  R.,  Constrained reaction coordinate dynamics for the simulation of
  rare events. \textit{Chem. Phys. Lett.} \textbf{156} (1989),
  472--477.

\bibitem{debocs98} Dellago, C.; Bolhuis, P. G.; Csajka, F.  S.;
  Chandler, D. Transition path sampling and the calculation of rate
  constants. \textit{J. Chem. Phys.} \textbf{ 108} (1998), 1964--1977.

\bibitem{deboge01} Dellago, C.; Bolhuis, P. G.; Geissler, P. L.
  Transition Path Sampling, Submitted to: \textit{ Adv. Chem. Phys.}
  (2001).

\bibitem{ERV1} E, W.; Ren, W.; Vanden-Eijnden, E. String method
  for the study of rare events. \textit{Phys. Rev.  B}, in press
  (2001).

\bibitem{ERV2} E, W.; Ren, W.; Vanden-Eijnden, E. Energy Landscape
  and Thermally Activated Switching of Submicron-sized Ferromagnetic
  Elements. \textit{J. App. Phys.}, submitted.

\bibitem{ERV3} E, W.; Ren, W.; Vanden-Eijnden, E.  Probing
  Multi-Scale Energy Landscapes Using the String Method. \textit{Phys.
    Rev. Lett.}, submitted.

\bibitem{ERV4} E, W.; Ren, W.; Vanden-Eijnden, E. Transition
  pathways in complex systems. Theory and numerical methods. In
  preparation.

\bibitem{FW} Freidlin, M. I.; Wentzell, A. D. \textit{Random Perturbations
  of Dynamical Systems}, 2nd ed. Springer, 1998.

\bibitem{hatabo90} H\"anggi, P.; Talkner, P.; Borkovec, M.
  Reaction-rate theory: fifty years after Kramers. \textit{Rev. Mod.
    Physics} \textbf{62} (1990), 251--341.

\bibitem{jomijo98} J\'onsson, H.; Mills, G.; Jacobsen, K. W. Nudged
  elastic band method for finding minimum energy paths of transitions.
  In: Classical and Quantum Dynamics in Condensed Phase Simulations.
  Edited by: Berne, B. J.; Cicotti, G.; Coker, D. F. World Scientific,
  1998.

\bibitem{olel96} Olender, R.; Elber, R. Calculation of classical
  trajectories with a very large time step: Formalism and numerical
  examples. \textit{J. Chem. Phys.} \textbf{105} (1996), 9299--9315.

\bibitem{Ren} Ren, W. \textit{Numerical Methods for the Study of Energy
Landscapes and Rare Events}, thesis, New York University, 2002.

\end{thebibliography}
\end{document}